\documentclass{article}

\usepackage{amsmath,amssymb,latexsym,theorem,bbm}

\newcommand{\TS}{\textstyle}
\newcommand{\SC}{\scriptstyle}

\newcommand{\CC}{\mathbb{C}}
\newcommand{\DD}{\mathbb{D}}
\newcommand{\NN}{\mathbb{N}}
\newcommand{\RR}{\mathbb{R}}
\newcommand{\ZZ}{\mathbb{Z}}

\newcommand{\dd}{\mathrm{d}}
\newcommand{\slu}{{\SC\mathrm{lu}}}

\newcommand{\cA}{\mathcal{A}}
\newcommand{\cB}{\mathcal{B}}
\newcommand{\cC}{\mathcal{C}}

\newcommand{\csL}{{\SC\cL}}
\newcommand{\cF}{\mathcal{F}}
\newcommand{\cL}{\mathcal{L}}

\newcommand{\cU}{\mathcal{U}}
\newcommand{\cV}{\mathcal{V}}
\newcommand{\cW}{\mathcal{W}}

\newcommand{\cY}{\mathcal{Y}}

\newcommand{\tC}{\widetilde{C}}

\newcommand{\tcC}{\widetilde{\cC}}

\newcommand{\EE}{\mathsf{E}}

\newcommand{\sPP}{{\SC\mathsf{P}}}
\newcommand{\PP}{\mathsf{P}}

\newcommand{\var}{\operatorname{var}}
\newcommand{\Var}{\operatorname{Var}}

\newcommand{\tr}{\operatorname{tr}}

\newcommand{\nt}{\lfloor nt\rfloor}
\newcommand{\nT}{\lfloor nT\rfloor}

\newcommand{\stoch}{\stackrel{\sPP}{\longrightarrow}}
\newcommand{\distr}{\stackrel{\csL}{\longrightarrow}}
\newcommand{\lu}{\stackrel{\slu}{\longrightarrow}}

\newcommand{\bmid}{\,\big|\,}
\newcommand{\Bmid}{\Big|}
\newcommand{\BMID}{\,\bigg|\,}
\newcommand{\bone}{\mathbbm{1}}
\newcommand{\vare}{\varepsilon}

\renewcommand{\leq}{\leqslant}
\renewcommand{\geq}{\geqslant}

\numberwithin{equation}{section}

\theoremstyle{plain}
\newtheorem{Lem}{Lemma}[section]
\newtheorem{Thm}[Lem]{Theorem}

\newtheorem{Cor}[Lem]{Corollary}
\theorembodyfont{\rm}

\newtheorem{Ex}[Lem]{Example}

\newcommand{\proofend}{\hfill$\square$}

\begin{document}

\begin{center}
 {\bfseries\Large A note on weak convergence of random step processes} \\[5mm]

 {\sc\large M\'arton $\text{Isp\'any}^*$} {\large and}
 {\sc\large Gyula $\text{Pap}^{*,\diamond}$} \\[5mm]
* University of Debrecen,
 Faculty of Informatics, Pf.~12,
 H--4010 Debrecen, Hungary \\[2mm]
 e--mails: Ispany.Marton@inf.unideb.hu (M. Isp\'any),
           Pap.Gyula@inf.unideb.hu (G. Pap). \\[2mm]
$\diamond$ Corresponding author. \\[5mm]
{\sl April 14, 2009.}
\end{center}

\vskip0.2cm

\renewcommand{\thefootnote}{}

\vspace*{-10mm}

\begin{abstract}
First, sufficient conditions are given for a triangular array of random vectors
 such that the sequence of related random step functions converges towards a
 (not necessarily time homogeneous) diffusion process.
These conditions are weaker and easier to check than the existing ones in the
 literature, and they are derived from a very general semimartingale
 convergence theorem due to Jacod and Shiryaev, which is hard to use directly.

Next, sufficient conditions are given for convergence of stochastic integrals
 of random step functions, where the integrands are functionals of the
 integrators.
This result covers situations which can not be handled by existing ones.
\end{abstract}

\section{Introduction}

The aim of the present paper is to obtain a useful theorem concerning
 convergence of step processes towards a diffusion process.
We derive sufficient conditions (see Theorem \ref{Conv2Diff} and Corollary
 \ref{Conv2DiffCor}) from a very general semimartingale convergence theorem
 due to Jacod and Shiryaev \cite[Theorem IX.3.39]{JSH}. 
(This theorem of Jacod and Shiryaev is hard to use directly, since one has to
 check the local strong majoration hypothesis, the local condition on big
 jumps, local uniqueness for the associated martingale problem, and the
 continuity condition.)
Theorem \ref{Conv2Diff} can also be considered as a generalization of the
 sufficient part of the functional martingale central limit theorem (see, e.g.,
 Jacod and Shiryaev \cite[Theorem VII.3.4]{JSH}), but Theorem \ref{Conv2Diff}
 allows not necessarily time homogeneous diffusion limit processes as well.
Similarly, Corollary \ref{Conv2DiffCor} can be considered as a generalization
 of the sufficient part of the Lindeberg-Feller functional central limit
 theorem (see, e.g., Jacod and Shiryaev \cite[Theorem VII.5.4]{JSH}).

There are several diffusion approximations in the literature, but they contain
 assumptions which are stronger and more complicated to check. 
For example, Ethier and Kurtz \cite[Theorem 7.4.1]{EK} deals only with the
 time homogeneous case, and their conditions (4.3)---(4.7) are hard to check.
The result of Joffe and M\'etivier \cite[Theorem 3.3.1]{JM} is not easy to use,
 since their conditions $(\mathrm{H}_1)$ and $(\mathrm{H}_4)$ are rather
 complicated to check.
Gikhman and Skorokhod \cite[Theorem 9.4.1]{GS} covers only convergence of
 Markov chains, and it contains Lipschitz conditions on the drift and diffusion
 coefficient of the limiting diffusion process, and assumes finite \ $2+\delta$
 \ moments for some \ $\delta>0$.
\ Our Theorem \ref{Conv2Diff} and Corollary \ref{Conv2DiffCor} are valid not
 only for martingales or Markov chains, since we do not suppose any dependence
 structure. 
The conditions are natural, since uniform convergence on compacts in
 probability (ucp) is involved.
(The role of the topology of the ucp is nicely explained by Kurtz and Protter
 \cite{KP3}.)

We also develope sufficient conditions (see Theorem \ref{Conv2Int} and
 Corollary \ref{Conv2IntCor}) for convergence of stochastic integrals of random
 step functions, where the integrand is a functional of the integrator.
We mention that our result covers situations which can not be handled by the
 convergence theorems of Jacod and Shiryaev \cite[Theorem IX.5.12,
 Theorem IX.5.16, Corollary IX.5.18, Remark IX.5.19]{JSH}.
There is a nice theory of convergence of stochastic integrals due to
 Jakubowski, M\'emin  and Pag\`es \cite{JMP} and to Kurtz and Protter
 \cite{KP1}, \cite{KP2}, \cite{KP3}.
The key result of this theory says that if \ $(\cU^n)_{n\in\NN}$ \ is a
 uniformly tight sequence of semimartingales (or, equivalently, it has
 uniformly controlled variations) then it is good in the sense that
 \ $(\cU^n,\cV^n,\cY^n)\distr(\cU,\cV,\cY)$ \ whenever
 \ $(\cU^n,\cV^n)\distr(\cU,\cV)$, \ where
 \ $\cY_t^n:=\int_0^t\cV_{s-}^n\,\dd\cU_s^n$ \ and
 \ $\cY_t:=\int_0^t\cV_{s-}\,\dd\cU_s$.
\ In our Theorem \ref{Conv2Int} and Corollary \ref{Conv2IntCor}, the sequence
 \ $(\cU^n)_{n\in\NN}$ \ of semimartingales is not necessarily good (see
 Example \ref{Ex_UT_UCV}).

In the proofs the simple structure of the approximating step processes and the
 almost sure continuity of the limiting diffusion process play a crucial role. 

As an application of these results, a Feller type diffusion approximation can
 be derived for critical multitype branching processes with immigration if the
 offspring mean matrix is primitive, and the asymptotic behavior of the
 conditional least squares estimator of the offspring mean matrix may be
 established, see Isp\'any and Pap \cite{IP}, which will be the content of a
 forthcoming paper.

\section{Convergence of step processes to diffusion processes}

A process \ $(\cU_t)_{t\in\RR_+}$ \ with values in \ $\RR^d$ \ is called a
 diffusion process if it is a weak solution of a stochastic differential
 equation
 \begin{equation}\label{SDE}
  \dd \, \cU_t 
  = \beta(t,\cU_t) \, \dd t + \gamma (t,\cU_t) \, \dd \cW_t,
  \qquad t\in\RR_+,
 \end{equation}
 where \ $\RR_+$ \ denotes the set of nonnegative real numbers,
 \ $\beta:\RR_+\times\RR^d\to\RR^d$ \ and
 \ $\gamma:\RR_+\times\RR^d\to\RR^{d\times r}$ \ are Borel functions and
 \ $(\cW_t)_{t\in\RR_+}$ \ is an \ $r$-dimensional standard Wiener process.

If \ $(\Omega,\cA,\PP)$ \ is a probability space, \ $\cF\subset\cA$ \ is a
 \ $\sigma$-algebra, and \ $\xi:\Omega\to\RR^d$ \ is a random variable with
 \ $\EE(\|\xi\|^2 \mid \cF)<\infty$ \ then \ $\Var(\xi\mid\cF)$ \ will denote
 the conditional variance matrix defined by
 \[
   \Var( \xi \mid \cF )
   := \EE \Big( \big( \xi - \EE( \xi \mid \cF ) \big)
                \big( \xi - \EE( \xi \mid \cF ) \big)^\top \bmid \cF \Big).
 \]
(Here and in the sequel, \ $\|x\|$ \ denotes the Euclidean norm of a (column)
 vector \ $x\in\RR^d$, \ $A^\top$ \ and \ $\tr A$ \ denote the transpose and
 the trace of a matrix \ $A$, \ respectively.)
The set of all nonnegative integers and the set of all positive integers will
 be denoted by \ $\ZZ_+$ \ and \ $\NN$, \ respectively.
The lower integer part and the positive part of \ $x\in\RR$ \ will be denoted
 by \ $\lfloor x\rfloor$ \ and \ $x_+$, \ respectively.

\begin{Thm}\label{Conv2Diff}
Let \ $\beta:\RR_+\times\RR^d\to\RR^d$ \ and
 \ $\gamma:\RR_+\times\RR^d\to\RR^{d\times r}$ \ be continuous functions.
Assume that the SDE \eqref{SDE} has a unique weak solution with \ $\cU_0=u_0$
 \ for all \ $u_0\in\RR^d$.
\ Let \ $\eta$ \ be a probability measure on \ $\RR^d$, \ and let
 \ $(\cU_t)_{t\in\RR_+}$ \ be a solution of \eqref{SDE} with initial
 distribution \ $\eta$.
\ For each \ $n\in\NN$, \ let \ $(U^n_k)_{k\in\ZZ_+}$ \ be a sequence of random
 variables with values in \ $\RR^d$ \ adapted to a filtration
 \ $(\cF^n_k)_{k\in\ZZ_+}$.
\ Let
 \[
   \cU^n_t := \sum_{k=0}^{\nt} U^n_k\,, \qquad t\in\RR_+, \quad n\in\NN.
 \]
Let \ $h:\RR^d\to\RR^d$ \ be a continuous function with compact support
 satisfying \ $h(x)=x$ \ in a neighborhood of 0.
Suppose \ $U^n_0\distr\eta$, \ and for each \ $T>0$,
 \begin{enumerate}
  \item[\textup{(i)}]
        $\sup\limits_{t\in[0,T]}
          \left\| \sum\limits_{k=1}^{\nt}
                   \EE ( h(U^n_k) \mid \cF^n_{k-1} )
                  - \int_0^t \beta(s,\cU^n_s) \, \dd s \right\|
         \stoch 0$,\\
  \item[\textup{(ii)}]
        $\sup\limits_{t\in[0,T]}
         \left\| \sum\limits_{k=1}^{\nt}
                  \Var ( h(U^n_k) \mid \cF^n_{k-1} )
                 - \int_0^t
                    \gamma(s,\cU^n_s) \gamma(s,\cU^n_s)^\top
                    \dd s \right\|
         \stoch 0$,\\
  \item[\textup{(iii)}]
        $\sum\limits_{k=1}^{\lfloor nT \rfloor}
          \PP ( \|U^n_k\| > \theta \mid \cF^n_{k-1} )
         \stoch 0$
        \ for all \ $\theta>0$.
 \end{enumerate}
Then \ $\cU^n \distr \cU$ \ as \ $n\to\infty$, \ i.e., the distributions of
 \ $\cU^n$ \ on the Skorokhod space \ $\DD(\RR^d)$ \ converge weakly to the
 distribution of \ $\cU$ \ on \ $\DD(\RR^d)$.
\end{Thm}

\noindent
\textbf{Proof.}
The process \ $(\cU_t)_{t\in\RR_+}$ \ is a semimartingale with characteristics
 \ $(\cB,\cC,0)$, \ where \ $\cB_t:=\int_0^t\beta(s,\cU_s)\,\dd s$, 
 \ $\cC_t:=\int_0^t\gamma(s,\cU_s)\gamma(s,\cU_s)^\top\dd s$ \ (see Jacod and
 Shiryaev \cite[III. \S\,2c]{JSH}).
In general, \ $\var\cB$ \ and \ $\tr\cC$ \ do not necessarily satisfy 
 majoration hypothesis, where \ $\var\alpha$ \ denotes the total variation of a
 function \ $\alpha\in\DD(\RR^d)$.
\ So we fix \ $T>0$, \ and stop the characteristics at \ $T$, \ that is, we
 consider the processes \ $\big(\cB^T_t\big)_{t\in\RR_+}$ \ and
 \ $\big(\cC^T_t\big)_{t\in\RR_+}$ \ defined by
 \[
   \cB^T_t := \int_0^{t\wedge T} \beta(s,\cU_s) \, \dd s, \quad
   \cC^T_t := \int_0^{t\wedge T} \gamma(s,\cU_s) \gamma(s,\cU_s)^\top \dd s,
 \]
 where \ $t\wedge T:=\inf\{t,T\}$.
\ Clearly, the stopped process \ $\big(\cU^T_t\big)_{t\in\RR_+}$ \ defined by
 \ $\cU^T_t:=\cU_{t\wedge T}$ \ is a semimartingale with characteristics
 \ $\big(\cB^T,\cC^T,0\big)$.

We will also consider the stopped processes
 \ $\big(\cU^{n,T}_t\big)_{t\in\RR_+}$, \ $n\in\NN$, \ defined by
 \ $\cU^{n,T}_t:=\cU^n_{t\wedge T}$.
\ We will check that all hypotheses of Theorem IX.3.39 of Jacod and Shiryaev
 \cite{JSH} are fulfilled.

Firstly, we check the local strong majoration hypothesis.
For each \ $a>0$, \ consider the mapping \ $\tau_a:\DD(\RR^d)\to[0,\infty]$
 \ defined by 
 \ $\tau_a(\alpha)
    :=\inf\{t\in\RR_+:\text{$|\alpha(t)|\geq a$ or $|\alpha(t-)|\geq a$}\}$
 \ for \ $\alpha\in\DD(\RR^d)$, \ where \ $\inf\emptyset:=\infty$.
\ Then the stopped processes
 \ $\big(\var\cB^T_{t\wedge\tau_a(\cU^T)}\big)_{t\in\RR_+}$ \ and
 \ $\big(\tr\cC^T_{t\wedge\tau_a(\cU^T)}\big)_{t\in\RR_+}$ \ are strongly
 majorized by the functions \ $t\mapsto b_{a,T}t$ \ and \ $t\mapsto c_{a,T}t$
 \ respectively, where
 \[
   b_{a,T} := \sup_{t\in[0,T]} \, \sup_{\|x\|\leq a}
               \|\beta(t,x)\|, \qquad
   c_{a,T} := \sup_{t\in[0,T]} \, \sup_{\|x\|\leq a}
               \|\gamma(t,x)\|^2.
 \]
Indeed, for all \ $s,t\in\RR_+$ \ with \ $s<t$, \ we have
 \begin{align*}
  \var\cB^T_{t\wedge\tau_a(\cU^T)}
  -\var\cB^T_{s\wedge\tau_a(\cU^T)}
  &=\int_{s\wedge T\wedge\tau_a(\cU^T)}
           ^{t\wedge T\wedge\tau_a(\cU^T)}
        \|\beta(u,\cU^T_u)\|\,\dd u,\\
  \tr\cC^T_{t\wedge\tau_a(\cU^T)}
  -\tr\cC^T_{s\wedge\tau_a(\cU^T)}
  &=\int_{s\wedge T\wedge\tau_a(\cU^T)}
        ^{t\wedge T\wedge\tau_a(\cU^T)}
     \|\gamma(u,\cU^T_u)\|^2\,\dd u.
 \end{align*}
The process \ $\big(\cU^T_t\big)_{t\in\RR_+}$ \ is a.s.\ continuous, hence
 \ $u\leq t\wedge T\wedge\tau_a(\cU^T)$ \ implies \ $\|\cU^T_u\|\leq a$
 \ a.s, thus
 \[
   \|\beta(u,\cU^T_u)\|\leq b_{a,T} \quad \text{a.s.}, \qquad
   \|\gamma(u,\cU^T_u)\|^2\leq c_{a,T} \quad \text{a.s.}
 \]
Consequently
 \begin{align*}
  \int_{s\wedge T\wedge\tau_a(\cU^T)}^{t\wedge T\wedge\tau_a(\cU^T)}
   \|\beta(u,\cU^T_u)\|\,\dd u
  & \leq b_{a,T}t-b_{a,T}s \qquad\text{a.s.,} \\
  \int_{s\wedge T\wedge\tau_a(\cU^T)}^{t\wedge T\wedge\tau_a(\cU^T)}
    \|\gamma(u,\cU^T_u)\|^2\,\dd u
  & \leq c_{a,T}t-c_{a,T}s \qquad\text{a.s.,}
 \end{align*} 
 hence the local strong majoration hypothesis holds.

The local condition on big jumps is obviously satisfied, since the third
 characteristic of the semimartingale \ $\big(\cU^T_t\big)_{t\in\RR_+}$ \ is
 \ 0.
By the assumption, the martingale problem associated to the characteristics
 \ $(\cB^T,\cC^T,0)$ \ admits a unique solution for each initial value
 \ $u_0\in\RR^d$, \ thus Theorem III.2.40 of Jacod and Shiryaev \cite{JSH}
 yields local uniqueness for the corresponding martingale problem as in
 Corollary III.2.41.
The continuity conditions are clearly implied by the continuity of the
 functions \ $\beta$ \ and \ $\gamma$.
\ Convergence of the initial distributions holds trivially.

For each \ $n\in\NN$, \ the stopped process
 \ $\big(\cU^{n,T}_t\big)_{t\in\RR_+}$ \ is also a semimartingale with
 characteristics
 \begin{align*}
  \cB^{n,T}_t & := \sum_{k=1}^{\lfloor n(t\wedge T)\rfloor}
                    \EE ( h(U^n_k) \mid \cF^n_{k-1} ), \\
  \cC^{n,T}_t & := 0, \\
  \nu^{n,T}([0,t]\times g)
  & := \sum_{k=1}^{\lfloor n(t\wedge T)\rfloor}
        \EE \big( g(U^n_k) \bone_{\{U^n_k\not=0\}} \bmid \cF^n_{k-1} \big)
 \end{align*}
 for \ $g:\RR^d\to\RR_+$ \ Borel functions, and modified second characteristic
 \[
   \tcC^{n,T}_t := \sum_{k=1}^{\lfloor n(t\wedge T)\rfloor}
                    \Var ( h(U^n_k) \mid \cF^n_{k-1} )
 \]
 (see Jacod and Shiryaev \cite[II.3.14, II.3.18]{JSH}).
For each \ $a>0$, \ assumptions (i)---(iii) imply  
 \begin{align*} 
  & \sup_{t\in[0,T]}
     \left\| \cB^{n,T}_{t\wedge\tau_a(\cU^T)}
             - \int_0^{t\wedge\tau_a(\cU^T)} \beta(s,\cU^n_s) \, \dd s\right\|
    \stoch 0, \\
  & \sup_{t\in[0,T]}
     \left\| \tcC^{n,T}_{t\wedge\tau_a(\cU^T)}
             - \int_0^{t\wedge\tau_a(\cU^T)}
                \gamma(s,\cU^n_s) \gamma(s,\cU^n_s)^\top \dd s\right\|
    \stoch 0, \\
  & \nu^{n,T}\big([0,\tau_a(\cU^T)] \times g_c\big)
    \stoch 0 \qquad
    \text{for all \ $c>0$,} 
 \end{align*}
 where \ $g_c:\RR^d\to\RR_+$ \ is defined by
 \begin{equation}\label{gc}
  g_c(x) := (c\|x\|-1)_+ \wedge 1.
 \end{equation}
(Indeed, \ $g_c(x) \leq \bone_{\{ \|x\| > 1/c \}}$ \ for all \ $x\in\RR^d$).
\ Therefore all hypotheses of Theorem IX.3.39 of Jacod and Shiryaev \cite{JSH}
 are fulfilled, hence for all \ $T>0$, \ $\cU^{n,T}\distr\cU^T$.
\ This implies that the finite dimensional distributions of the processes
 \ $\cU^n$ \ converge to the corresponding finite dimensional distributions of
 the process \ $\cU$ \ (see Jacod and Shiryaev \cite[VI.3.14]{JSH}).

The aim of the following discussion is to show the tightness of
 \ $\{\cU^n:n\in\NN\}$, \ which will imply \ $\cU^n\distr\cU$.
\ For each \ $T>0$, \ by Prokhorov's Theorem, convergence
 \ $\cU^{n,T}\distr\cU^T$ \ implies tightness of \ $\{\cU^{n,T}:n\in\NN\}$.
\ By Theorem VI.3.21 of Jacod and Shiryaev \cite{JSH}, this implies
 \begin{align*} 
  & \PP \Big( \sup_{t\in[0,T]} \| \cU^{n,T}_t \| > K \Big) \to 0
    \qquad \text{as \ $n\to\infty$ \ and \ $K\to\infty$,} \\
  & \PP \left( w_T' \big( \cU^{n,T}, \theta \big) > \delta \right) \to 0
    \qquad
    \text{as \ $n\to\infty$ \ and \ $\theta \downarrow 0$ \ for all
          \ $\delta>0$,}
 \end{align*}
 where \ $w_T'(\alpha,\cdot)$ \ denotes the ``modulus of continuity'' on
 \ $[0,T]$ \ for a function \ $\alpha\in\DD(\RR^d)$ \ (see Jacod and Shiryaev
 \cite[VI.1.8]{JSH}).
Since the above convergences hold for all \ $T>0$, \ we conclude for all
 \ $T>0$ \ that
 \begin{align*} 
  & \PP \Big( \sup_{t\in[0,T]} \| \cU^n_t \| > K \Big) \to 0
    \qquad \text{as \ $n\to\infty$ \ and \ $K\to\infty$,} \\
  & \PP \left( w_T' \big( \cU^n, \theta \big) > \delta \right) \to 0
    \qquad
    \text{as \ $n\to\infty$ \ and \ $\theta \downarrow 0$ \ for all
          \ $\delta>0$.}
 \end{align*} 
Again by Theorem VI.3.21 of Jacod and Shiryaev \cite{JSH}, this implies
 tightness of \ $\{\cU^n:n\in\NN\}$, \ and we obtain \ $\cU^n\distr\cU$.
\proofend

\begin{Cor}\label{Conv2DiffCor}
Let \ $\beta$, \ $\gamma$, \ $\eta$, \ $(U^n_k)_{k\in\ZZ_+}$,
 \ $(\cF^n_k)_{k\in\ZZ_+}$ \ and \ $\cU^n$ \ for \ $n\in\NN$ \ be as in
 Theorem \ref{Conv2Diff}.
Suppose that \ $\EE \big( \|U^n_k\|^2 \bmid \cF^n_{k-1} \big) < \infty$ \ for
 all \ $n,k\in\NN$.
\ Assume that the SDE \eqref{SDE} has a unique weak solution with \ $\cU_0=u_0$
 \ for all \ $u_0\in\RR^d$.
\ Let \ $(\cU_t)_{t\in\RR_+}$ \ be a solution of \eqref{SDE} with initial
 distribution \ $\eta$.
\ Suppose \ $U^n_0\distr\eta$, \ and for each \ $T>0$,
 \begin{enumerate}
  \item [\textup{(i)}]
         $\sup\limits_{t\in[0,T]}
          \left\| \sum\limits_{k=1}^{\nt}
                   \EE ( U^n_k \mid \cF^n_{k-1} )
                  - \int_0^t \beta(s,\cU^n_s) \, \dd s\right\|
          \stoch 0$,\\
  \item [\textup{(ii)}]
        $\sup\limits_{t\in[0,T]}
         \left\| \sum\limits_{k=1}^{\nt}
                  \Var ( U^n_k \mid \cF^n_{k-1} )
                 - \int_0^t
                    \gamma(s,\cU^n_s) \gamma(s,\cU^n_s)^\top
                    \dd s \right\|
         \stoch 0$,\\
  \item [\textup{(iii)}]
        $\sum\limits_{k=1}^{\lfloor nT \rfloor}
          \EE \big( \|U^n_k\|^2 \bone_{\{\|U^n_k\| > \theta\}}
                    \bmid \cF^n_{k-1} \big)
         \stoch 0$
        \ for all \ $\theta>0$.
 \end{enumerate}
Then \ $\cU^n \distr \cU$ \ as \ $n\to\infty$.
\end{Cor}

\noindent
\textbf{Proof.}
Clearly, there exists \ $K\geq 1$ \ such that \ $h(x)=x$ \ for
 \ $\|x\|\leq1/K$, \ $h(x)=0$ \ for \ $\|x\|\geq K$, \ and \ $\|h(x)\|\leq K$
 \ for all \ $x\in\RR^d$.
\ Hence \ $h(x)-x=0$ \ for \ $\|x\|\leq1/K$ \ and
 \ $\|h(x)-x\| \leq \|h(x)\| + \|x\| \leq K + \|x\| \leq (K^2+1)\|x\|$ \ for 
 \ $\|x\|\geq1/K$.
\ Thus, we conclude 
 \begin{equation}\label{ineq}
  \|h(x)-x\|
  \leq (K^2+1) \|x\| \bone_{\{\|x\|\geq1/K\}}
  \leq (K^2+1)K \|x\|^2 \bone_{\{\|x\|\geq1/K\}}
 \end{equation}
 for all \ $x\in\RR^d$.
\ For all \ $T>0$ \ and all \ $t\in[0,T]$, \ applying \eqref{ineq}, we get 
 \begin{align*} 
  \left\| \sum_{k=1}^{\nt}
             \EE ( h(U^n_k) \mid \cF^n_{k-1} )
          - \sum_{k=1}^{\nt}
             \EE ( U^n_k \mid \cF^n_{k-1} ) \right\|
  &\leq \sum_{k=1}^{\nT}
           \EE \big( \|h(U^n_k)- U^n_k\| \bmid \cF^n_{k-1} \big) \\
  & \leq (K^2+1)K
         \sum_{k=1}^{\nT} 
          \EE \Big( \|U^n_k\|^2 \bone_{\{\|U^n_k\|\geq1/K\}} 
                    \mid \cF^n_{k-1} \Big),
 \end{align*}
 which together with assumptions (i) and (iii) of this corollary imply
 condition (i) of Theorem \ref{Conv2Diff}.
We have
 \begin{align*} 
  \Var \big( h(U^n_k) \bmid \cF^n_{k-1} \big)
   - \Var \big( U^n_k \bmid \cF^n_{k-1} \big)
   = \EE \big( h(U^n_k) h(U^n_k)^\top - U^n_k (U^n_k)^\top 
               \bmid \cF^n_{k-1} \big)\\
    + \left( \EE \big( h(U^n_k) \bmid \cF^n_{k-1} \big)
             \EE \big( h(U^n_k)^\top \bmid \cF^n_{k-1} \big)
             - \EE \big( U^n_k \bmid \cF^n_{k-1} \big)
               \EE \big( (U^n_k)^\top \bmid \cF^n_{k-1} \big)
      \right).
 \end{align*}
For arbitrary matrices \ $A,B,C,D\in\RR^{d\times r}$, \ we have
 \[
   \|AB^\top-CD^\top\|
   \leq \|A-C\| \cdot \|B\| + \|A\| \cdot \|B-D\| + \|A-C\| \cdot \|B-D\|,
 \]
 hence applying \eqref{ineq} and \ $\|h(x)\|\leq K$ \ valid for all
 \ $x\in\RR^d$, \ we obtain
 \begin{align*} 
  \sum_{k=1}^{\nt}
   \left\| \EE \big( h(U^n_k) h(U^n_k)^\top - U^n_k (U^n_k)^\top 
                     \bmid \cF^n_{k-1} \big)
   \right\|
  & \leq \sum_{k=1}^{\nt}
          \EE \big( 2 \| h(U^n_k) - U^n_k \| \| h(U^n_k) \| 
                    + \| h(U^n_k) - U^n_k \|^2 \bmid \cF^n_{k-1} \big)\\
  & \leq (K^2+1)(3K^2+1)
         \sum_{k=1}^{\nt}
          \EE \big( \| U^n_k \|^2 \bone_{\{\| U^n_k \| \geq 1/K \}}  
                    \bmid \cF^n_{k-1} \big).
 \end{align*}
In a similar way, we obtain
 \begin{align*} 
  &\sum_{k=1}^{\nt}
    \left\| \EE \big( h(U^n_k) \bmid \cF^n_{k-1} \big)
            \EE \big( h(U^n_k)^\top \bmid \cF^n_{k-1} \big)
            - \EE \big( U^n_k \bmid \cF^n_{k-1} \big)
              \EE \big( (U^n_k)^\top \bmid \cF^n_{k-1} \big)
    \right\| \\
  & \leq 2K^2(K^2+1)
         \sum_{k=1}^{\nt}
          \EE \big( \| U^n_k \|^2 \bone_{\{\| U^n_k \| \geq 1/K \}}  
                    \bmid \cF^n_{k-1} \big)
         + K^2(K^2+1)^2
           \left( \sum_{k=1}^{\nt}
                   \EE \big( \| U^n_k \|^2 \bone_{\{\| U^n_k \| \geq 1/K \}}  
                             \bmid \cF^n_{k-1} \big)
           \right)^2.
 \end{align*}
These inequalities together with assumptions (ii) and (iii) of this corollary
 imply condition (ii) of Theorem \ref{Conv2Diff}. 
We have
 \[
   \PP \big( \| U^n_k \| > \theta \bmid \cF^n_{k-1} \big)
   \leq \theta^{-2}
        \EE \Big( \|U^n_k\|^2 \bone_{\{\|U^n_k\|\geq\theta\}}
                  \mid \cF^n_{k-1} \Big),
 \]
 thus assumption (iii) of this corollary implies (iii) of Theorem
 \ref{Conv2Diff}.
\proofend

\begin{Ex}\label{Ex_UT_UCV}
We give an example for a system \ $(U^n_k)_{n\in\NN,\,k\in\ZZ_+}$ \ of random
 variables satisfying conditions (i)---(iii) of Corollary \ref{Conv2DiffCor},
 such that the sequence \ $(\cU^n)_{n\in\NN}$ \ of semimartingales is not good
 (see the Introduction).

Let \ $(\eta_k)_{k\in\NN}$ \ be independent standard normal random variables.
Let \ $U^n_0:=0$, \ $U^n_{3j}:=-\eta_j/\sqrt{n}$,
 \ $U^n_{3j-1}:=U^n_{3j-2}:=\eta_j/\sqrt{n}$ \ and
 \ $\cF^n_{j-1}:=\sigma(U^n_0,\dots,U^n_{j-1})$ \ for \ $j,n\in\NN$.
\ Then conditions (i)---(iii) of Corollary \ref{Conv2DiffCor} are satisfied
 with \ $\beta=0$ \ and \ $\gamma=1/\sqrt{3}$.
\ For each \ $n\in\NN$, \ let
 \[
   \int_0^t \cU^n_{s-} \,\dd\cU^n_s
   = \sum_{k=1}^{\nt} U^n_k \sum_{j=1}^{k-1} U^n_j
   = \frac{1}{2} (\cU^n_t)^2
     - \frac{1}{2} \sum_{k=1}^{\nt} (U^n_k)^2.
 \]
Then, by Corollary \ref{Conv2DiffCor}, \ $\cU^n\distr\cU:=\cW/\sqrt{3}$,
 \ where \ $(\cW_t)_{t\in\RR_+}$ \ is a standard Wiener process.
Moreover, \ $\sum_{k=1}^{\nt} (U^n_k)^2 \stoch t$, \ hence
 \ $\int_0^t \cU^n_{s-} \, \dd\cU^n_s 
    \distr \frac{1}{6} (\cW_t)^2 -\frac{1}{2}t$.
 \ But, by It\^o's formula,
 \ $\int_0^t \cU_{s-} \,\dd \cU_s = \frac{1}{6} ((\cW_t)^2-t)$,
 \ thus the sequence
 \ $\left(\int_0^t \cU^n_{s-} \,\dd\cU^n_s\right)_{n\in\NN}$
 \ does not converge to \ $\int_0^t \cU_{s-} \,\dd\cU_s$.
\ Consequently, the sequence \ $(\cU^n)_{n\in\NN}$ \ of semimartingales is not
 good.
\end{Ex}

\section{Convergence of integrals of step processes}

For a function \ $\alpha\in\DD(\RR^d)$ \ and for a sequence
 \ $(\alpha_n)_{n\in\NN}$ \ in \ $\DD(\RR^d)$, \ we write \ $\alpha_n\lu\alpha$
 \ if \ $(\alpha_n)_{n\in\NN}$ \ converges to \ $\alpha$ \ locally uniformly,
 i.e., if \ $\sup_{t\in[0,T]}\|\alpha_n(t)-\alpha(t)\|\to0$ \ as \ $n\to\infty$
 \ for all \ $T>0$.
\ The space of all continuous functions \ $\alpha:\RR_+\to\RR^d$ \ will be
 denoted by \ $\CC(\RR^d)$.
\ For measurable mappings \ $\Phi:\DD(\RR^d)\to\DD(\RR^p)$ \ and
 \ $\Phi_n:\DD(\RR^d)\to\DD(\RR^p)$, \ $n\in\NN$, \ we will denote by
 \ $C_{\Phi,(\Phi_n)_{n\in\NN}}$ \ the set of all functions
 \ $\alpha\in\CC(\RR^d)$ \ such that \ $\Phi(\alpha)\in\CC(\RR^p)$ \ and
 \ $\Phi_n(\alpha_n)\lu\Phi(\alpha)$ \ whenever \ $\alpha_n\lu\alpha$ \ with
 \ $\alpha_n\in\DD(\RR^d)$, \ $n\in\NN$.
\ If \ $\Phi_n=\Phi$ \ for all \ $n\in\NN$ \ then we write simply \ $C_\Phi$
 \ instead of \ $C_{\Phi,(\Phi_n)_{n\in\NN}}$.
\ Further, \ $\tC_{\Phi,(\Phi_n)_{n\in\NN}}$ \ will denote the set of all
 functions \ $\alpha\in C_{\Phi,(\Phi_n)_{n\in\NN}}$ \ such that
 \ $\Phi(\alpha_n)\lu\Phi(\alpha)$ \ whenever \ $\alpha_n\lu\alpha$ \ with
 \ $\alpha_n\in\DD(\RR^d)$, \ $n\in\NN$.
\ Finally, \ $D_{\Phi,(\Phi_n)_{n\in\NN}}$ \ will denote the set of all
 functions \ $\alpha\in\DD(\RR^d)$ \ such that
 \ $\Phi_n(\alpha_n)\to\Phi(\alpha)$ \ in \ $\DD(\RR^p)$ \ whenever
 \ $\alpha_n\to\alpha$ \ in \ $\DD(\RR^d)$ \ with \ $\alpha_n\in\DD(\RR^d)$,
 \ $n\in\NN$.
\ We need the following version of the continuous mapping theorem several
 times.
 
\begin{Lem}\label{Conv2Funct}
Let \ $(\cU_t)_{t\in\RR_+}$ \ and \ $(\cU^n_t)_{t\in\RR_+}$, \ $n\in\NN$, \ be
 stochastic processes with values in \ $\RR^d$ \ such that \ $\cU^n\distr\cU$.
\ Let \ $\Phi:\DD(\RR^d)\to\DD(\RR^p)$ \ and
 \ $\Phi_n:\DD(\RR^d)\to\DD(\RR^p)$, \ $n\in\NN$, \ be measurable mappings
 such that \ $\PP\big(\cU\in C_{\Phi,(\Phi_n)_{n\in\NN}}\big)=1$.
\ Then \ $\Phi_n(\cU^n) \distr \Phi(\cU)$.
\end{Lem}

\noindent
\textbf{Proof.}
In view of the continuous mapping theorem (see, e.g., Billingsley
 \cite[Theorem 5.5]{BI}), it suffices to check that
 \ $\PP\big(\cU\in D_{\Phi,(\Phi_n)_{n\in\NN}}\big)=1$.
\ For a function \ $\alpha\in\CC(\RR^d)$, \ $\alpha_n\to\alpha$ \ in
 \ $\DD(\RR^d)$ \ if and only if \ $\alpha_n\lu\alpha$ \ (see, e.g., Jacod and
 Shiryaev \cite[VI.1.17]{JSH}).
Consequently,
 \ $\CC(\RR^d)\cap\Phi^{-1}(\CC(\RR^p))\cap D_{\Phi,(\Phi_n)_{n\in\NN}}
    =C_{\Phi,(\Phi_n)_{n\in\NN}}$
 \ implying
 \ $D_{\Phi,(\Phi_n)_{n\in\NN}}\supset C_{\Phi,(\Phi_n)_{n\in\NN}}$.
\proofend

\begin{Thm}\label{Conv2Int}
Let \ $\beta$, \ $\gamma$, \ $\eta$, \ $(U^n_k)_{k\in\ZZ_+}$,
 \ $(\cF^n_k)_{k\in\ZZ_+}$ \ and \ $\cU^n$ \ for \ $n\in\NN$ \ be as in
 Theorem \ref{Conv2Diff}.
Assume that the SDE \eqref{SDE} has a unique weak solution with \ $\cU_0=u_0$
 \ for all \ $u_0\in\RR^d$.
\ Let \ $(\cU_t)_{t\in\RR_+}$ \ be a solution of \eqref{SDE} with initial
 distribution \ $\eta$.

For each \ $n\in\NN$ \ and \ $k\in\ZZ_+$, \ let
 \ $\psi_{n,k}:(\RR^d)^{k+1}\to\RR^p$ \ be a Borel function, and let
 \ $\Psi_n:\DD(\RR^d)\to\DD(\RR^p)$ \ be defined by
 \[
   \Psi_n(\alpha)(t)
   := \psi_{n,\nt}
       \Big( \alpha\big({\TS\frac{1}{n}}\big)-\alpha(0), \dots,
             \alpha\big({\TS\frac{\nt}{n}}\big)
              - \alpha\big({\TS\frac{\nt-1}{n}}\big) \Big)
 \]
 for \ $\alpha\in\DD(\RR^d)$.
\ Let
 \begin{align*}
  V^n_k & := \psi_{n,k} ( U^n_0, \dots, U^n_k ),
   \qquad k\in\ZZ_+, \quad n\in\NN, \\
  \cV^n_t& := V^n_{\nt} = \Psi_n(\cU^n)_t,
   \qquad t\in\RR_+, \quad n\in\NN, \\
  \cY^n_t& := \sum_{k=1}^{\nt} V^n_{k-1} \otimes U^n_k
            = \int_0^t \cV^n_{s-} \otimes \dd \, \cU^n_s,
   \qquad t\in\RR_+, \quad n\in\NN.
 \end{align*}
Let \ $\Psi:\DD(\RR^d)\to\DD(\RR^p)$ \ be a measurable mapping such that
 \ $\PP\big(\cU\in\tC_{\Psi,(\Psi_n)_{n\in\NN}}\big)=1$.
\ Let 
 \[ 
  \cV_t := \Psi(\cU)_t, \qquad
  \cY_t := \int_0^t \cV_{s-} \otimes \dd \, \cU_s,
  \qquad t\in\RR_+. 
 \]
Let the mappings \ $\beta':\DD(\RR^d)\to\DD(\RR^d\times\RR^{pd})$ \ and 
 \ $\gamma':\DD(\RR^d)\to\DD(\RR^{d\times r}\times\RR^{(pd)\times r})$ \ be
 defined by
 \[
  \beta'(\alpha)(s) := \begin{bmatrix}
                        \beta(s,\alpha(s)) \\
                        \Psi(\alpha)(s) \otimes \beta(s,\alpha(s))
                       \end{bmatrix},
  \quad
  \gamma'(\alpha)(s) := \begin{bmatrix}
                         \gamma(s,\alpha(s)) \\
                         \Psi(\alpha)(s) \otimes \gamma(s,\alpha(s))
                        \end{bmatrix}.
 \] 
Let \ $h':\RR^d\times\RR^{pd}\to\RR^d\times\RR^{pd}$ \ be a continuous function
 with compact support satisfying \ $h'(x)=x$ \ in a neighborhood of 0.
Suppose \ $U^n_0\distr\eta$, \ and for each \ $T>0$,
 \begin{enumerate}
  \item [\textup{(i)}]
        $\sup\limits_{t\in[0,T]}
          \left\| \sum\limits_{k=1}^{\nt}
                   \EE ( h'(U^n_k,V^n_{k-1} \otimes U^n_k) \mid \cF^n_{k-1} )
                  - \int_0^t \beta'(\cU^n)_s \, \dd s \right\|
         \stoch 0$,\\
  \item [\textup{(ii)}]
        $\sup\limits_{t\in[0,T]}
         \left\| \sum\limits_{k=1}^{\nt}\!
                  \Var ( h'(U^n_k,V^n_{k-1} \!\otimes U^n_k)
                         \!\mid\! \cF^n_{k-1} )
                 \!-\! \int_0^t\!
                    \gamma'(\cU^n)_s \gamma'(\cU^n)_s^\top
                    \dd s \right\|
         \stoch 0$,\\
  \item [\textup{(iii)}]
        $\sum\limits_{k=1}^{\lfloor nT \rfloor}
          \PP \big( \| U^n_k \| (1 + \| V^n_{k-1} \|) > \theta 
                    \bmid \cF^n_{k-1} \big)
         \stoch 0$
        \ for all \ $\theta>0$.
 \end{enumerate}
Then \ $(\cU^n,\cV^n,\cY^n) \distr (\cU,\cV,\cY)$ \ as \ $n\to\infty$.
\end{Thm}

\noindent
\textbf{Proof.}
Our first aim is to prove \ $(\cU^n,\cY^n)\distr(\cU,\cY)$.
\ We start by showing that the sequence \ $(\cU^n,\cY^n)_{n\in\NN}$ \ is tight
 in \ $\DD(\RR^d\times\RR^{pd})$, \ and for this we will use Theorem VI.4.18 of
 Jacod and Shiryaev \cite{JSH}.
By the assumptions, the sequence \ $(\cU^n_0,\cY^n_0)=(U^n_0,0)$, \ $n\in\NN$,
 \ is weakly convergent, thus obviously tight in \ $\RR^d\times\RR^{pd}$,
 \ hence condition (i) of Theorem VI.4.18 of Jacod and Shiryaev \cite{JSH}
 holds.
For each \ $n\in\NN$, \ the process \ $(\cU^n_t,\cY^n_t)_{t\in\RR_+}$ \ is a
 semimartingale with characteristics \ $(\cB'^n,\cC'^n,\nu'^n)$ \ relative to
 the truncation function \ $h'$ \ given by
 \begin{align*}
  \cB'^n_t & := \sum_{k=1}^{\nt}
                 \EE ( h'(U^n_k,V^n_{k-1} \otimes U^n_k)
                       \mid \cF^n_{k-1} ), \\
  \cC'^n_t & := 0, \\
  \nu'^n([0,t]\times g)
  & := \sum_{k=1}^{\nt}
        \EE \Big( g(U^n_k,V^n_{k-1} \otimes U^n_k)
                  \bone_{\{(U^n_k,V^n_{k-1} \otimes U^n_k)\not=0\}}
                  \mid \cF^n_{k-1} \Big)
 \end{align*}
 for \ $g:\RR^d\times\RR^{pd}\to\RR_+$ \ Borel functions, and modified
 second characteristic
 \[
   \tcC'^n_t := \sum_{k=1}^{\nt}
                 \Var ( h'(U^n_k,V^n_{k-1} \otimes U^n_k) \mid \cF^n_{k-1} )
 \]
 (see Jacod and Shiryaev \cite[II.3.14, II.3.18]{JSH}).
For all \ $T>0$, \ $\theta>0$, \ $\vare>0$,
 \begin{align*}
  &\PP\Big(\nu'^n\big([0,T]\times\bone_{\{\|x\|>\theta\}}\big)>\vare\Big)\\
  &=\PP\bigg(\sum_{k=1}^{\nT}
              \PP(\|(U^n_k,V^n_{k-1} \otimes U^n_k)\|>\theta
                  \mid\cF^n_{k-1})>\vare\bigg)
   \to0
 \end{align*}
 by assumption (iii), hence condition (ii) of Theorem VI.4.18 of Jacod and
 Shiryaev \cite{JSH} holds.

In order to check condition (iii) of Theorem VI.4.18 of Jacod and Shiryaev
 \cite{JSH}, first we will show
 \begin{equation}\label{ConvBn'}
  \int_0^t\beta'(\cU^n)_s\,\dd s\distr\int_0^t\beta'(\cU)_s\,\dd s
  \qquad \text{in \ $\DD(\RR^d\times\RR^{pd})$.}
 \end{equation}
We will apply Lemma \ref{Conv2Funct}.
We have \ $\int_0^t\beta'(\cU)_s\,\dd s=\Phi_{\beta'}(\cU)_t$, \ and for each
 \ $n\in\NN$, \ $\int_0^t\beta'(\cU^n)_s\,\dd s=\Phi_{\beta'}(\cU^n)_t$
 \ with the measurable mapping
 \ $\Phi_{\beta'}:\DD(\RR^d)\to\DD(\RR^d\times\RR^{pd})$ \ given by
 \[
   \Phi_{\beta'}(\alpha)(t) := \int_0^t\beta'(\alpha)(s)\,\dd s,
   \qquad \alpha\in\DD(\RR^d), \quad t\in\RR_+.
 \]
Observe that assumptions (i)--(iii) imply that conditions (i)--(iii) of
 Theorem \ref{Conv2Diff} hold, thus we conclude \ $\cU^n\distr\cU$ \ as
 \ $n\to\infty$.
\ In order to show \ $\PP\big(\cU\in C_{\Phi_{\beta'}}\big)=1$, \ it is enough
 to check \ $C_{\Phi_{\beta'}}\supset\tC_{\Psi,(\Psi_n)_{n\in\NN}}$. 
\ Clearly \ $\Phi_{\beta'}(\CC(\RR^d))\subset\CC(\RR^d\times\RR^{pd})$. 
\ Now we fix \ $T>0$, \ a function \ $\alpha\in\tC_{\Psi,(\Psi_n)_{n\in\NN}}$
 \ and a sequence \ $(\alpha_n)_{n\in\NN}$ \ in \ $\DD(\RR^d)$ \ with
 \ $\alpha_n\lu\alpha$.
\ Obviously
 \[
   \sup_{t\in[0,T]} \|\Phi_{\beta'}(\alpha_n)-\Phi_{\beta'}(\alpha)\|
   \leq T \sup_{t\in[0,T]} \|\beta'(\alpha_n)(t)-\beta'(\alpha)(t)\|,
 \]
 hence it suffices to show
 \begin{gather}
  \sup_{t\in[0,T]} \|\beta(t,\alpha_n(t))-\beta(t,\alpha(t))\| \to 0,
   \label{Unif1} \\
  \sup_{t\in[0,T]} \|\Psi(\alpha_n)(t)\otimes\beta(t,\alpha_n(t))
                     -\Psi(\alpha)(t)\otimes\beta(t,\alpha(t))\| \to 0.
   \label{Unif2}
 \end{gather}
For sufficiently large \ $n\in\NN$, \ we have
 \ $\sup_{t\in[0,T]} \|\alpha_n(t)-\alpha(t)\| \leq 1$, \ thus
 \ $\sup_{t\in[0,T]} \|\alpha_n(t)\|
    \leq 1 + \sup_{t\in[0,T]} \|\alpha(t)\|<\infty$.
\ The function \ $\beta$ \ is uniformly continuous on the compact set
 \ $[0,T]\times\{x\in\RR^d:\|x\|\leq 1 + \sup_{t\in[0,T]} \|\alpha(t)\|\}$,
 \ hence \eqref{Unif1} holds.
Moreover,
 \begin{align*}
  &\|\Psi(\alpha_n)(t)\otimes\beta(t,\alpha_n(t))
     -\Psi(\alpha)(t)\otimes\beta(t,\alpha(t))\| \\
  &\leq \|\Psi(\alpha_n)(t)-\Psi(\alpha)(t)\| \|\beta(t,\alpha_n(t))\|
        +\|\beta(t,\alpha_n(t))-\beta(t,\alpha(t))\| \|\Psi(\alpha)(t)\|.
 \end{align*}
Continuity of \ $\Psi(\alpha)$ \ implies
 \ $\sup_{t\in[0,T]}\|\Psi(\alpha)(t)\|<\infty$.
\ For sufficiently large \ $n\in\NN$,
 \ $\sup_{t\in[0,T]}\|\beta(t,\alpha_n(t))\|
    \leq 1 + \sup_{t\in[0,T]}\|\beta(t,\alpha(t))\| < \infty$
 \ (by convergence \eqref{Unif1} and by continuity of \ $\alpha$ \ and
 \ $\beta$).
By \ $\Psi(\alpha_n)\lu\Psi(\alpha)$, \ \eqref{Unif2} is also satisfied, and we
 conclude \ $C_{\Phi_{\beta'}}\supset\tC_{\Psi,(\Psi_n)_{n\in\NN}}$.
\ Consequently, \ $\PP\big(\cU\in C_{\Phi_{\beta'}}\big)=1$, \ and by Lemma
 \ref{Conv2Funct}, we obtain \eqref{ConvBn'}.
If \ $\alpha\in\CC(\RR^d)$ \ and \ $(\alpha_n)_{n\in\NN}$ \ is a sequence in
 \ $\DD(\RR^d)$ \ with \ $\alpha_n\lu\alpha$ \ then for all \ $T>0$,
 \ $\sup_{t\in[0,T]}\|\alpha_n(t)\|\to\sup_{t\in[0,T]}\|\alpha(t)\|$ \ as
 \ $n\to\infty$.
\ (See, e.g., Proposition VI.2.4 of Jacod and Shiryaev \cite{JSH}.)
Hence, by the continuous mapping theorem, we obtain
 \[
  \sup\limits_{t\in[0,T]}
   \left\| \int_0^t \beta'(\cU^n)_s \, \dd s
           - \int_0^t \beta'(\cU)_s \, \dd s \right\|
   \distr 0
   \qquad \text{as \ $n\to\infty$.}
 \]
This together with assumption (i) implies
 \begin{equation}\label{SupConvBn'}
  \sup\limits_{t\in[0,T]}
   \left\| \cB'^n_t - \int_0^t \beta'(\cU)_s \, \dd s \right\|
   \stoch 0
   \qquad \text{as \ $n\to\infty$ \ for all \ $T>0$.}
 \end{equation}
Particularly, the sequence \ $(\cB'^n)_{n\in\NN}$ \ is \ $C$-tight in
 \ $\DD(\RR^d\times\RR^{pd})$.
\ Indeed, the Skorokhod topology is coarser than the local uniform topology,
 hence \eqref{SupConvBn'} implies
 \ $\varrho\big(\cB'^n,\Psi_{\beta'}(\cU)\big)\stoch0$, \ where \ $\varrho$
 \ denotes a distance on \ $\DD(\RR^d)$ \ compatible with the Skorokhod
 topology.
Consequently, \ $\cB'^n\distr\Psi_{\beta'}(\cU)$ \ with
 \ $\PP(\Psi_{\beta'}(\cU)\in\CC(\RR^d\times\RR^{pd}))=1$.
\ In a similar way, the sequence \ $(\tcC'^n)_{n\in\NN}$ \ is \ $C$-tight in
 \ $\DD(\RR^{d\times r}\times\RR^{(pd)\times r})$.
\ Moreover, assumption (iii) yields
 \begin{equation}\label{SupConvNu'}
  \nu'^n([0,T]\times g_c)\stoch0
  \qquad \text{as \ $n\to\infty$}
 \end{equation}
 for all \ $T>0$ \ and all \ $c>0$, \ where \ $g_c:\RR^d\times\RR^{pd}\to\RR_+$
 \ is defined by \eqref{gc}.
Therefore all hypotheses of Theorem VI.4.18 of Jacod and Shiryaev \cite{JSH}
 are fulfilled, hence we conclude that the sequence \ $(\cU^n,\cY^n)_{n\in\NN}$
 \ is tight in \ $\DD(\RR^d\times\RR^{pd})$.

It remains to prove that if a sub-sequence, still denoted by
 \ $(\cU^n,\cY^n)_{n\in\NN}$, \ weakly converges to a limit distribution then
 the limit is the distribution of \ $(\cU,\cY)$.
\ For this we will apply Theorem IX.2.22 of Jacod and Shiryaev \cite{JSH}.
The process \ $(\cU_t,\cY_t)_{t\in\RR_+}$ \ is a semimartingale with
 characteristics \ $(\cB',\cC',0)$, \ where
 \[
  \cB'_t := \int_0^t \beta'(\cU)_s \,\dd s, 
  \qquad 
  \cC'_t := \int_0^t \gamma'(\cU)_s \gamma'(\cU)_s^\top \dd s
 \] 
(see Jacod and Shiryaev \cite[IX.5.3]{JSH}).
By Remark IX.2.23 of Jacod and Shiryaev \cite{JSH}, assumptions (i)--(iii) of
 Theorem \ref{Conv2Int} imply that condition (i) of Theorem IX.2.22 in
 \cite{JSH} is met.
To prove the continuity condition (ii) of Theorem IX.2.22 in \cite{JSH},
 consider the measurable mapping
 \ $\Phi:\DD(\RR^d)
         \to\DD\big(\RR^d
                    \times(\RR^d\times\RR^{pd})
                    \times(\RR^{d\times r}\times\RR^{(pd)\times r})\big)$
 \ given by
 \[
  \Phi(\alpha)(t)
  := \big(\alpha(t),\Phi_{\beta'}(\alpha)(t),\Phi_{\gamma'}(\alpha)(t)\big),
  \qquad \alpha\in\DD(\RR^d), \quad t\in\RR_+.
 \]
As we have already proved,
 \ $\PP\big(\cU\in C_{\Phi_{\beta'}}\cap C_{\Phi_{\gamma'}}\big)=1$.
\ The local uniform topology on \ $\DD(\RR^m)$ \ is the \ $m$-fold product of
 the local uniform topology on \ $\DD(\RR)$, \ hence we obtain
 \ $C_\Phi\supset C_{\Phi_{\beta'}}\cap C_{\Phi_{\gamma'}}$.
\ Using again that the Skorokhod topology is coarser than the local uniform
 topology, we conclude \ $D_\Phi\supset C_\Phi$.
\ Consequently, the continuity condition \ $\PP\big(\cU\in D_\Phi\big)=1$
 \ holds.
Hence all hypotheses of Theorem IX.2.22 of Jacod and Shiryaev \cite{JSH} are
 met, therefore \ $(\cU^n,\cY^n)\distr(\cU,\cY)$.
\ Again by Lemma \eqref{Conv2Funct}, we obtain
 \ $(\cU^n,\cV^n,\cY^n) \distr (\cU,\cV,\cY)$.
\proofend

\begin{Cor}\label{Conv2IntCor}
Let \ $\beta$, \ $\gamma$, \ $\eta$, \ $(U^n_k)_{k\in\ZZ_+}$,
 \ $(\cF^n_k)_{k\in\ZZ_+}$ \ and \ $\cU^n$ \ for \ $n\in\NN$ \ be as in Theorem
 \ref{Conv2Diff}.
Suppose that \ $\EE \big( \|U^n_k\|^2 \bmid \cF^n_{k-1} \big) < \infty$ \ for
 all \ $n,k\in\NN$.
\ Assume that the SDE \eqref{SDE} has a unique weak solution with
 \ $\cU_0=u_0$ \ for all \ $u_0\in\RR^d$.
\ Let \ $(\cU_t)_{t\in\RR_+}$ \ be a solution with initial distribution
 \ $\eta$.
\ Let \ $\Psi$, \ $\cV$, \ $\cY$, \ $\beta'$, \ $\gamma'$,
 \ $(\psi_{n,k})_{k\in\NN}$, \ $\Psi_n$, \ $(V^n_k)_{k\in\ZZ_+}$, \ $\cV^n$
 \ and \ $\cY^n$ \ for \ $n\in\NN$ \ be as in Theorem \ref{Conv2Int}.
Suppose that \ $\PP\big(\cU\in\tC_{\Psi,(\Psi_n)_{n\in\NN}}\big)=1$.
\ Suppose \ $U^n_0\distr\eta$, \ and for each \ $T>0$,
 \begin{enumerate}
  \item [\textup{(i)}]
         $\sup\limits_{t\in[0,T]}
          \left\| \sum\limits_{k=1}^{\nt}
                   \EE \left( \begin{bmatrix}
                               U^n_k \\
                               V^n_{k-1} \otimes U^n_k)
                              \end{bmatrix} 
                              \BMID \cF^n_{k-1} \right)
                  - \int_0^t
                     \beta'(\cU^n)_s \, \dd s\right\|
          \stoch 0$,\\
  \item [\textup{(ii)}]
        $\sup\limits_{t\in[0,T]}
         \left\| \sum\limits_{k=1}^{\nt}
                  \Var \! \left( \begin{bmatrix}
                                  U^n_k \\
                                  V^n_{k-1} \otimes U^n_k)
                                 \end{bmatrix}
                                 \BMID \cF^n_{k-1} \right)
                 - \int_0^t
                    \gamma'(\cU^n)_s \gamma'(\cU^n)_s^\top
                    \dd s \right\|
         \stoch 0$,\\
  \item [\textup{(iii)}]
        $\sum\limits_{k=1}^{\lfloor nT \rfloor}\!\!
          \EE \Big( \! \| U^n_k \|^2 (1 \!+\! \| V^n_{k-1} \|^2)
                       \bone_{\{\| U^n_k \| (1+ \| V^n_{k-1} \|)
                                 > \theta\}}
                       \Bmid \cF^n_{k-1} \! \Big)
         \!\stoch\! 0$ \ for all \ $\theta>0$.
 \end{enumerate}
Then \ $(\cU^n,\cV^n,\cY^n) \distr (\cU,\cV,\cY)$.
\end{Cor}

\noindent
\textbf{Proof.}
This follows from Theorem \ref{Conv2Int} in the same way as Corollary
 \ref{Conv2DiffCor} from Theorem \ref{Conv2Diff}.
\proofend

\noindent{\bf Acknowledgements.}
The authors have been supported by the Hungarian Scientific Research Fund
 under Grant No.\ OTKA T-048544 and OTKA T-079128.

\end{document}